\documentclass[amsart, 11pt]{article}

\usepackage[utf8]{inputenc}
\usepackage{cite}
\usepackage{makeidx}
\usepackage{latexsym}
\usepackage{amsmath}
\usepackage{amssymb}
\usepackage{amscd}
\usepackage{bbm}
\usepackage{mathrsfs}
\usepackage{mathtools}
\usepackage{theorem}
\usepackage[all]{xy}

\DeclareMathAlphabet{\mathbbold}{U}{bbold}{m}{n}

\newtheorem{teo}{Theorem}[section]
\newtheorem{propo}[teo]{Proposition}
\newtheorem{lema}[teo]{Lemma}
\newtheorem{coro}[teo]{Corollary}

{\theorembodyfont{\rmfamily }
\newtheorem{defi}[teo]{Definition}

\newtheorem{obser}[teo]{Remark}
\newtheorem{ejem}[teo]{Example}
}
\def\demo{\noindent \textit{Proof: }}

\def\qed{\hspace*{\fill }$\square $}

\def\ch{\mathrm {ch}}

 \DeclareMathOperator{\rk}{rank}
 
\DeclareMathOperator{\Td}{Td}

\def\11{1}
\def\11{\mathbbm{1}}

\def\cf{\emph{cf. }}

\def\pt{*}

\def\raya{\ \underline{\phantom{a}}\ }

\def\qed{\hspace*{\fill }$\square $ }

\def\BB{\mathsf{B}}

\def\O{\mathcal{O}}
 \def\PP{\mathbb{P}}
                    \def\QQ{\mathbb{Q}}
                     
                                                            \def\SS{\mathrm{S}}
 \def\TT{\mathsf{T}}

\begin{document}

\title{\sc The Riemann-Roch theorem for the Adams operations}
\author{A. Navarro \thanks{a.navarrogarmendia@icmat.es,  ICMAT, Spain}, J. Navarro \thanks{navarrogarmendia@unex.es, Universidad de Extremadura, Spain \newline Supported in part by    H2020-MSCA-IF-897784,  ERDF and Junta de Extremadura - Consejer\'{\i}a de Econom\'{\i}a, Ciencia y Agenda Digital programs IB18087 and GR21055.} }

\maketitle

\begin{abstract}
We prove the classical Riemann-Roch theorems for the Adams operations $\,\psi^j\,$ on $K$-theory: a statement with coefficients on $\mathbb{Z}[j^{-1}]$,  that holds for arbitrary projective morphisms,  as well as another one with integral coefficients,  that is valid for closed immersions.  In presence of rational coefficients, we also analyze the relation between the corresponding Riemann-Roch formula for one Adams operation and the analogous formula for the Chern character.

To do so,  we complete the elementary exposition of the work of Panin-Smirnov 
that was initiated by the first author in a previous work.  Their notion of oriented cohomology theory of algebraic varieties allows to use classical arguments to prove general and neat statements, which imply all the aforementioned results as particular cases.
\end{abstract}



\bigskip

On a compact Riemann surface, a pioneering work of B. Riemann (\cite{Riemann}) proved the lower bound
$$ l(D) \, \geq \, {\rm deg}\, D - g + 1 $$
for the dimension $l(D)$ of the complex vector space of meromorphic functions with zeroes and poles prescribed by a divisor  $D$. 
This inequality unveiled a strong relation between such a complex analytic quantity, $l(D)$, and purely topological invariants; namely, the degree of the divisor and the topological genus $g$ of the surface. Not long after that,  Riemann's student G. Roch identified (\cite{Roch}) the difference between both sides of the inequality: if $K$ denotes the canonical divisor on the surface, 
$$l(D) - l(K-D) \, =  \, {\rm deg} D - g + 1 \ . $$

After a series of generalizations during the first half of the twentieth century,  A. Grothendieck's approach (\cite{SGA6},\cite{BorelSerre}) entailed a Copernican revolution on the understanding of that formula.  First of all, he introduced the $K$-group $K(X)$ of vector bundles of an algebraic $k$-variety $X$ and realized that it has a cohomological behaviour: it is a contravariant functor of rings and, for any projective morphism $f \colon Y \to X$ between smooth varieties, there exits an exceptional direct image $f_! \colon K(Y) \to K(X)$. 
Also, $K$-theory is related to the Chow ring with rational coefficients via a morphisms of functor of rings $\mathrm{ch}\colon K \to CH^\bullet_\mathbb{Q}$, called the Chern character.
The generalization of the Riemann-Roch formula says how to modify this Chern character in order to make it compatible with the corresponding direct images: the correction term is the Todd series of the virtual tangent bundle $T_f := T_Y - f^* T_X \in K(Y)$. 

To be precise,  the Grothendieck-Riemann-Roch theorem \cite{BorelSerre} establishes that for any proper morphism $f \colon Y \to X $ between smooth quasi-projective $k$-varieties, the following square commutes:
\begin{equation}\label{CuadradoGRR}
\vcenter{\xymatrix{
K(Y)  \ar[rrr]^-{f_!} \ar[d]_{ \ch} & & & K(X) \ar[d]^{\ch} \\
CH^\bullet (Y)_\mathbb{Q} \ar[rrr]^-{f_* (\Td (T_f) \cdot \raya  ) } & & & CH^\bullet (X)_\mathbb{Q}}}
\end{equation}

On the other hand, $K$-theory comes equipped with certain natural transformations of rings, $\psi^j \colon K \to K$, called the Adams operations. 
If we replace the Chern character by these Adams operations, there exists an analogous formula (\cite[Theorem 7.6]{FultonLang},  \cite[Theorem 16.6]{Manin}, \cite{PinkRossler}): for any proper morphism $f \colon Y \to X $ between smooth quasi-projective $k$-varieties, the following square commutes:

\begin{equation}\label{CuadradoARR}
\vcenter{\xymatrix{
K(Y)  \ar[rrr]^-{f_!} \ar[d]_{\psi^j} & & & K(X) \ar[d]^{\psi^j} \\
K(Y) [j^{-1}]  \ar[rrr]^-{f_! ( \theta^j (\Omega_f)^{-1} \cdot \raya ) } & & & K(X) [j^{-1}]}} 
\end{equation}
where $K(X) [j^{-1}] := K(X) \otimes_\mathbb{Z} \mathbb{Z}[j^{-1}]$, $\Omega_f = T_f^* \in K(Y) $ is the relative cotangent bundle and $\,\theta^j \colon K \to K [j^{-1}]^*\,$ is Bott's cannibalistic class. 


The purpose of this work is to expound how the ideas  that have emerged in motivic homotopy theory (mainly, in the works of I. Panin, A. Smirnov and J. Riou,  see  \cite{Panin-Smirnov}, \cite{Panin}, \cite{Riou}, \cite{Smirnov}) can shed light into 
the theorem above, even without any use of the motivic framework.  In this sense, this paper continues and completes — and may well be read in conjunction with -- the exposition in  \cite{ExpositionesI}. We also hope that our paper will make this result more accessible for the reader not willing to invest in the formalism developed in the standard reference \cite{FultonLang}.  Hence, we claim no originality in the following results, apart from the exposition. 

As a brief account of the relevance of this Adams-Riemann-Roch formula, let us mention that both squares, (\ref{CuadradoGRR}) and (\ref{CuadradoARR}),  when written with rational coefficients,  were proved to be ``equivalent'' in an abstract setting (\cite[Chapter III, \S 4]{FultonLang}).  Inspired by that,  C.  Soul\'{e} used the Adams operations in higher $K$-theory to propose the first definition of motivic cohomology with rational coefficients (\cite{Soule}). 

Let us now comment on the plan of this article. First of all, we consider oriented cohomology theories over algebraic varieties,  in the sense of Panin-Smirnov.  This notion 
is wide enough to include at the same time examples such as all Weyl cohomologies, the Chow ring, Grothendieck's $K$-theory of vector bundles or Levine-Morel's algebraic cobordism. 

In the second place, any morphism of functor of rings $\varphi \colon A \to \bar{A}$ between oriented cohomology theories has an associated formal power series $\SS $ with coefficients on the cohomology ring of the one point space $\bar{A}(\pt)$ (\cf Definition \ref{SerieAsociada}).  This series encodes the behaviour of $\varphi$ on the cohomology of projective spaces; in the case of the Chern character $\mathrm{ch}\colon K(X)\to CH^\bullet (X)_\mathbb{Q}$,  this series is the inverse of the usual Todd series (Example \ref{EjemSeries}.2). We then show how the standard arguments commonly used to prove Riemann-Roch formulas (separated study of regular immersions and projections, deformation to the normal bundle...) in fact apply to arbitrary morphism of rings between arbitrary cohomology theories, establishing:

\medskip
{\bf Theorem \ref{RR}:} 
{\it Let $\,\varphi \colon A\to \bar{A}\,$ be a morphism of functor of rings between oriented cohomology theories, whose associated series $\SS \in \bar{A}(\pt)[[t]] $ is invertible.  For any projective morphism $\,f \colon Y \to X\,$ between smooth, quasi-projective $k$-varieties, the following diagram commutes
\begin{equation*}
\vcenter{\xymatrix{
A(Y)  \ar[rrr]^-{f_*} \ar[d]_{\varphi} & & & A(X) \ar[d]^{\varphi} \\
\bar{A}(Y) \ar[rrr]^-{f_* ( \SS_{\times}^{-1} (T_f) \cdot  \raya )  } & & & \bar{A} (X)}} 
\end{equation*} where $\SS_\times^{-1} \colon K \to \bar{A}$ denotes the multiplicative extension of the series $\SS^{-1}$.} 
\medskip



Then, we explain how this theorem subsumes both the Grothendieck-Riemann-Roch and the Adams-Riemann-Roch formulae presented before (Corollary \ref{RRA} and Corollary \ref{RRG}).

Also, as an application of this level of generality,  we prove Corollary \ref{TeoremaCuboMagico}, a statement for which we know no reference,  and that may be understood as another incarnation of the  ``equivalence",  in presence of rational coefficients,  between the squares (\ref{CuadradoGRR}) and (\ref{CuadradoARR})  (see \cite[Ch. III, \S 4] {FultonLang} or \cite[App. B]{Duma}).

Finally, it is also worth pointing out that,  prior to all that, we separately analyze the case of regular closed immersions, where we prove this finer statement:

\medskip
{\bf Theorem \ref{Main}:}
{\it Let $\,\varphi \colon A \to \bar{A}\,$ be a morphism of functor of rings between oriented cohomology theories, with associated series $\, \SS \in \bar{A}(\pt)[[t]]\,$.  For any closed inmersion $\, i \colon Z \hookrightarrow X\,$, the following square commutes:
\begin{equation*}
\vcenter{\xymatrix{
A(Z)  \ar[rrr]^-{i_*} \ar[d]_{\varphi} & & & A(X) \ar[d]^{\varphi} \\
\bar{A} (Z) \ar[rrr]^-{i_* (\SS_{\times} (N_i) \cdot  \raya ) } & & & \bar{A}(X)}}
\end{equation*}
where $\, N_i := i^* T_{X} - T_Z = - T_i\,$ denotes the normal bundle of $i$. }
\medskip

\noindent This result produces, as particular cases, the Adams-Riemann-Roch formula without denominators of \cite[Ch. V, \S 6]{FultonLang}  (Corollary \ref{CoroRRAsindeno}),
as well as an integral formula for a twisted Chern character that we have not found in the literature (Corollary \ref{FORMULA_NUEVA}).




\section{Oriented cohomology theories on algebraic varieties}

We work in the setting of oriented cohomology theories on algebraic varieties, in the sense of Panin-Smirnov (\cite{Panin-Smirnov}). Let us briefly recall this formalism, as exposed in \cite{ExpositionesI}.

\begin{defi}\label{DefiCohomologia} Let $k$ be a field.  An \textit{oriented cohomology theory} is a contravariant functor $A\,$ from the category of smooth quasi-projective varieties over $k\,$ into the category of commutative rings. Part of the data are also \textit{direct images for projective morphisms}; that is to say,  a morphism of groups $f_*\colon A(Y)\to A(X)$, for any projective morphism $f\colon Y\to X$, such that
$\mathrm{Id}_*=\mathrm{Id}$, $(fg)_*=f_*g_*\,$ and such that the projection formula
$f_*(f^*(x)y)=xf_*(y)$ holds.

On an oriented cohomology theory $\,A$, any smooth closed subvariety $\,i\colon Z\to X\,$ has a \textit{fundamental class}
$$\,[Z]^A\coloneqq i_*(1)\in A(X)\ ,$$ and any line
bundle $\,L\to X\,$ has a \textit{Chern class} $$\,c_1^A(L)\coloneqq s_0^*(s_{0*}(1))\in A(X) \ , $$ where $\,s_0\colon X\to L\,$ is the zero section.

These data are also assumed to satisfy the following axioms:
\begin{enumerate}
\item
The ring morphism $i_1^*+i_2^*\colon A(X_1\amalg X_2)\to A(X_1)\oplus A(X_2)$ is an
isomorphism, where $i_j\colon X_j\to X_1\amalg X_2$ is the natural immersion (in particular, $A(\emptyset )=0$).

\item For any affine bundle
$\pi \colon P\to X$, 
the ring morphism $\pi ^*\colon A(X)\to A(P)$ is an isomorphism.

\item
For any smooth closed subvariety $i\colon Z \to X$,  the following is an exact sequence: $$\
A(Z)\xrightarrow{i_*}A(X)\xrightarrow{j^*}A(X-Z)\ .$$

\item
If a morphism $f\colon \bar X\to X$ is transversal to a smooth closed subvariety
$i\colon Z\to X$ of codimension $d$, 
then the following square is commutative:
$$
\xymatrix{A(Z) \ar[r]^-{f^*} \ar[d]^{i _*} & A(\bar Z) \ar[d]^{i _*}\\ A(X)
\ar[r]^-{f^*} & A(\bar X). }
$$

\item
For any morphism $f\colon
Y\to X$ and any projective bundle $\pi \colon \mathbb{P}(E)\to X$,  the following square is commutative:
$$
\xymatrix{A(\mathbb{P}(E))\ar[r]^-{f^*} \ar[d]^{\pi _*} & A(\mathbb{P}(f^*E))
\ar[d]^{\pi _*}\\ A(X) \ar[r]^-{f^*} & A(Y). }
$$

\item If $E\to X$ is a vector bundle of rank $r+1$, then the $A(X)$-module in $A(\mathbb{P}(E))$ is free of rank $r+1$:
$$
A(\mathbb{P}(E))=A(X)\oplus A(X)x_E\oplus \ldots \oplus A(X)x_E^r.
$$ where $x_E=c_1^A(\xi
_E)$ is the first Chern class of the tautological
line bundle.
\end{enumerate}


A \textit{morphism of cohomologies} $\varphi \colon A\to \bar A\,$ is a natural transformation preserving direct images; 
i.e.,  it is a morphism of functors of rings $\,\varphi \colon A\to \bar A\,$ such that
$$
\varphi (f_*(a))=\bar f_*(\varphi (a))
$$
for any projective morphism $f\colon Y\to X$ and any element $a\in A(Y)$.
\end{defi}

\medskip

\begin{ejem}
\begin{enumerate}
\item
Weil cohomologies are oriented cohomologies in the sense above. These examples include singular cohomology of complex varieties, de Rham cohomology of varieties over a field $k$ of characteristic zero, or $l$-adic cohomology, with $k$ of characteristic different from $l\,$.

\item
The Chow ring $\, CH^\bullet (X)\,$ of rational equivalence classes of algebraic cycles is also an oriented cohomology theory (\cite{Fulton}).

\item
Grothendieck's $K$-theory is an oriented cohomology theory where the cohomology rings $\,K(X)\,$ have no graduation (v. gr., \cite{ExpositionesI}).


\item
Algebraic cobordism $\, \Omega^{\bullet} (X)$, as constructed by Levine and Morel (\cite{LevineMorel}), is also an oriented cohomology.
\end{enumerate}
\end{ejem} 

The following fact, commonly known as {\it Jouanolou's Trick}, assures that any quasi-projective $k$-variety is ``homotopically equivalent" to an affine variety:

\begin{teo}[\cite{Jouanolou}]\label{Jouanolou}
For any quasi-projective $k$-variety $\,X$, there exists an affine bundle $\,P \to X\,$ such that $\,P\,$ is an affine variety. 
\end{teo}

\noindent In particular, for any given line bundle $L \to X$, there exists an affine bundle $\,p \colon P \to X\,$ and a morphism $\,\varphi \colon P \to \mathbb{P}^d\,$ such that $\, p^* L = \varphi^* (\xi^*_d)$, where $\, \xi_d \to \mathbb{P}^d\,$ stands for the tautological line bundle. 


\begin{defi}\label{FGL}
A {\it commutative formal group law} over a ring $\,R\,$ is a power series $\,f(u,v) \in R[[u,v]]\,$ that satisfies the axioms for a commutative group operation with $0$ as the neutral element; namely,
\begin{itemize}
\item[a)] $f(u,v) = f(v,u) \in R[[u,v]]$. 
\item[b)] $f(u,0) = u \in R[[u]] $.
\item[c)] $f(u, f(v,w)) = f(f(u,v), w) \in R[[u,v,w]]$. 
\item[d)] There exist a series $\, opp(u) \in R[[u]]\,$ such that $\, opp(0) = 0\,$ and $\,f(u, opp(u) ) = 0 \in R[[u]]$. 
\end{itemize}
The {\it additive group law} is $f(u , v)= u + v\ ,$ and the {\it multiplicative} 
{\it group law} is $f(u ,v)= u+v-uv$.
\end{defi}


\medskip
On any oriented cohomology theory $A$,  the cohomology rings of projective spaces satisfy the projective bundle theorem (\cite[Corollary 1.3]{ExpositionesI}); that is to say, the maps $c_1^A(\xi^*_d) \mapsto u$   define isomorphisms of rings, for all $d \geq 1$, 
$$A(\PP^d)=A(\pt)[u]/(u^{d+1}) \ .$$ 

As a consequence, 

\begin{propo}\label{LeyGrupoFormal}
Let $\,A\,$ be an oriented cohomology theory.  There exists a unique commutative formal group law $\,f^A(u,v)\in \,A(\pt)[[u,v]]\,$ such that
$$ c_1^A (L_x \otimes L_y) = f^A(x,y) $$ for any pair of line bundles $L_x$, $L_y$, with Chern classes $x := c_1^A (L_x)$ and $y := c_1^A (L_y)$.

We say that $A$ {\it follows the law} $f^A$, or that the series $f^A(u,v)$ is the {\it formal group law associated to $A$}. 
\end{propo}

\demo Axiom 6 of Definition \ref{DefiCohomologia} implies that the map $  c_1 ( \pi^*_1 \xi^*_d ) \mapsto u$, $ c_1 (\pi_2^* \xi^*_d) \mapsto v $ is an isomorphism of rings: 
\[ A \left(\mathbb{P}^d \times \mathbb{P}^d \right) = 
A(\pt) [u,v] / (u^{d+1} , v^{d+1}) \ . \]

These isomorphisms are compatible with the restrictions $i^* \colon \mathbb{P}^d \to \mathbb{P}^{d-1}$; therefore, there exists a unique series $f (u,v) \in A(\pt) [[u,v]]$ such that, for any $d$,
$$ c_1 ( \pi_1^* \xi^*_d \otimes \pi_2^* \xi^*_d ) = f ( u, v) \ . $$

Now, for any pair of line bundles $L_x$ and $L_y$ on any variety $X$,  use Theorem  \ref{Jouanolou} to assure the existence of an affine bundle $p\colon P \to X$ and maps $h_1, h_2 \colon P \to \mathbb{P}^d$ such that $p^* L_x = h_1^* \xi^*_d$ and $p^* L_y = h_2^* \xi^*_d$. 

Considering the map $h := h_1 \times h_2 \colon P \to \mathbb{P}^d \times \mathbb{P}^d$, we  can write
\begin{align*}
p^* \left( c_1 (L_x \otimes L_y ) \right) &= c_1 (p^* L_x \otimes p^*L_y ) =  c_1 ( h^* \pi_1^* \xi^*_d \otimes h^* \pi_2^* \xi^*_d ) = h^* \left( c_1 ( \pi_1^* \xi^*_d \otimes \pi_2^* \xi^*_d ) \right)  \\
&= h^* \left( f(u,v) \right) = f \left( h^* u , h^* v \right) = f(p^* x, p^* y) = p^* (f(x,y) )\ ,
\end{align*} and the series $f(x,y)$ satisfies the required property, because $p^*$ is an isomorphism.

The first three axioms of formal group law follow from the fact that $L_x \otimes L_y = L_y \otimes L_x$, $L_x \otimes \mathcal{O}_X = L_x$ and $ (L_x \otimes L_y) \otimes L_z = L_x \otimes (L_y \otimes L_z) $. 

As for the last one, let $opp(u)\in A(\pt)[[u]] $ be the unique series such that
\[ c_1( \xi_d) = opp(u) \in A(\pt)[u]/(u^{d+1}) \simeq A(\PP^d)\ , \qquad \forall \ d \geq 1 \ , \]
via the identification  $u\mapsto c_1 (\xi^*_d) $. The relation $\xi^*_d \otimes \xi_d = \mathcal{O}_{\mathbb{P}^d}$ implies that this series also satisfies
\[ f ( u, opp(u)) = 0 \in A(\pt) [u]/(u^{d+1}) \ , \qquad \forall \ d \geq 1 \ , \]
and hence $f (  u, opp(u)) = 0$ in $A(\pt) [[u]]$. 

\qed

\begin{ejem}\label{RelationXeY}

\begin{enumerate} 
\item
Graded cohomology theories follow the additive group law (\cite[Lemma 3.1]{ExpositionesI}). These examples include all Weil cohomologies or the Chow ring.

\item Grothendieck's $K$-theory follows the multiplicative group law,  whereas algebraic cobordism follows the universal group law (in particular,  $\, \Omega^{\bullet} (\pt) \,$ is the Lazard ring,  see \cite{LevineMorel}).

\item
The equality $\xi_d \otimes \xi^*_d = \mathcal{O}_X$ implies a relation between the Chern classes $x_d = c_1^A (\xi_d)$ and $y_d = c_1^A(\xi^*_d)$: if $A$ follows the additive law, then $x_d=- y_d$, whereas, if $A$ follows the multiplicative law, then $x_d =-y_d/(1-y_d)$.

\end{enumerate}
\end{ejem}

Morphisms of oriented cohomologies preserve Chern classes, so that:

\begin{coro}\label{CorolarioSerie}
If there exists a morphism of cohomologies $\, \varphi \colon A \to \bar{A}$, then both cohomology theories follow the same group law.
\end{coro}

\demo Let $f$ denote the formal group law associated to $A$. Then
\begin{equation*}
c_1^{\bar{A}}(L_x \otimes L_y ) \,  = \, \varphi (c_1^A(L_x \otimes L_y )) \, = \, \varphi (f(x,y)) \, = \, f(\varphi (x) , \varphi (y) ) \, = \, f( \bar{x} , \bar{y})  \ , 
\end{equation*}
where we write $\bar{x} = c_1^{\bar{A}} (L_x)$, and $\bar{y} = c_1^{\bar{A}} (L_y)$.

\qed

\medskip

In order to fix notation, let us recall some results of the theory of Chern classes for these cohomology theories: 

\begin{defi}
The {\it Chern classes} $c_n^A(E)\in A(X)$ of a vector bundle $E\to X$ of rank $r$ are the coefficients of the characteristic polynomial
$c(E)=x^r-c_1^A(E)x^{r-1}+\ldots +(-1)^rc_r^A(E)$ of the endomorphism of the free
$A(X)$-module $A(\mathbb{P}(E))$ defined by the multiplication by $x_E=c_1^A(\xi
_E)$, so that
$$
x_E^r-c_1^A(E)x_E^{r-1}+\ldots +(-1)^rc_r^A(E)=0 \ . 
$$
\end{defi}


These Chern classes are compatible with inverse images, and the total Chern class $\,c(E) := 1+c_1^A(E)t+c_2^A(E)t^2+\ldots
+c_r^A(E)t^r\,$ is additive (\cite[Theorem 1.1 and 1.2]{ExpositionesI}. Hence, it extends to a morphisms of functor of groups
$$ c \colon K \longrightarrow 1 + A\, [[t]] \ , $$ whose $n^{th}$-degree part is called the {\it $n^{th}$-Chern class map} $c_n\colon K \to A$.

The Chern classes so defined are nilpotent elements of the cohomology (\cite[Corollary 1.4]{ExpositionesI}). 

\section{Grothendieck's $K$-theory}\label{EjemploK}   

The main example for our interests is the oriented cohomology theory defined by the $K$-group of vector bundles (\cite{SGA6},  \cite{Weibel} or \cite{ExpositionesI}).  

In this example,  the fundamental class of a smooth closed subvariety $i\colon Y\to X$ is its structural sheaf:
$$
[Y]^K=i_*(1)=\mathcal{O}_Y\in K(X) \ , 
$$
and the Chern class of a line bundle $L \to X$ is the element $1-L^* \in K(X) $. 

Line bundles generate this $K$-group, in the following sense (\cite[Remark 2]{ExpositionesI}):

\medskip
\noindent {\bf Splitting Principle:} {\em Let $A$ be an oriented cohomology.
For any vector bundle $E\to X$, there exists a base change $\pi \colon X'\to X$ such that $\pi ^*\colon A(X)\to A(X')$ is
injective and $\pi ^*E=L_{\alpha _1}+\ldots +L_{\alpha _r}$ is a sum of line bundles
in $K(X')$.  }
\medskip

\noindent These cohomology classes $\, \alpha _1,\ldots ,\alpha _r \in A (X')\,$ appearing in the statement above are called the {\it roots} of $E$. 

\medskip

\subsubsection*{Morphism of groups $K \to A$}

Let $\,K^+(X)\,$ stand for the commutative monoid generated by the isomorphism classes of vector bundles over $\,X\,$, modulo the relations $\,[E] = [E'] + [E'']\,$ produced by exact sequences $\,0 \to E' \to E \to E'' \to 0\,$.

\begin{lema}\label{UnicidadMorfismoGrupos}
Let $A$ be an oriented cohomology. 
Morphisms of functor of monoids $\,\phi \colon K^+ \to (A , \cdot )$, as well as morphisms of functor of groups $\,\phi \colon K \to A\,$, are determined by their values on the tautological line bundles $\,\{\xi_d\}_{d\in \mathbb{N}}$.
\end{lema}

\demo Let us only detail the case of morphisms of groups. Let $\, \phi, \bar{\phi} \colon K \to A\,$ be morphisms of functor of groups that coincide on each $\xi_d$. 

By Theorem \ref{Jouanolou}, they coincide on any line bundle $\,L\to X$: let $p\colon P\to X$ be an affine bundle and $\,\varphi \colon P\to \mathbb{P}^d$ a morphism such that  $p^*L = \varphi^*(\xi_d)$. Then,
$$
p^* \phi (L) = \phi (p^* L) = \phi ( \varphi^* \xi_d) = \varphi^* \phi (\xi_d) = \varphi^* \bar{\phi} (\xi_d) = \bar{\phi} ( \varphi^* \xi_d) = \bar{\phi} ( p^* L) = p^* \bar{\phi}(L) \ , 
$$
and therefore $\,\phi (L)  = \bar{\phi} (L)$, because $p^* $ is injective.

As both $\phi$ and $ \bar{\phi}$ are functorial and additive, the Splitting Principle implies that they coincide on any vector bundle $E$ and, hence, on any element of the $K$-group.

\qed

\begin{propo}\label{ExtensionesDeSeries}
Let $A$ be an oriented cohomology and $\,F \in A(\pt)[[t]]\,$ be any series.

\begin{itemize}
\item {\bf Multiplicative extension}: There exists a unique morphism of functor of monoids, 
$$ 
F_\times \colon K^+ \longrightarrow (A , \cdot )  \ , 
$$ 
whose value on the the tautological line bundles is  $ F_\times (\xi_d) := F(x_d)$.

\item {\bf Additive extension}: There exists a unique morphism of functor of groups 
$$ F_+ \colon K \longrightarrow A   \ , $$ whose value on the tautological line bundles  is $\, F_+ (\xi_d)  := F(x_d)$.
\end{itemize} 
\end{propo}

\demo Uniqueness is a consequence of Lemma \ref{UnicidadMorfismoGrupos} above. To prove the existence, let us define, on a vector bundle $\,E\to X\,$ with roots $\alpha_1, \ldots , \alpha_r \in A(X')\,$,
\begin{align*}
F_{\times} (E) &:= F(\alpha_1) \cdot \ldots \cdot F(\alpha_r)  \\
F_{+} (E) &:= F(\alpha_1) + \ldots + F(\alpha_r)  
\end{align*}
where these elements of $A(X)$ are considered as power series in the elementary symmetric functions of the roots, that are precisely the Chern classes of $E$. 

These formulae are easily checked to be compatible with inverse images; $F_{\times}$ is a morphism of monoids and $F_{+}$ defines an additive function, so it extends to the whole $K$-group.

\qed

\begin{teo}\label{MorfismosDeGrupos}
Any morphism of functor of monoids $\, \phi \colon K^+ \to (A,\cdot) \,$ is the multiplicative extension, $\phi = F_\times$, of a unique series $\,F \in A(\pt)[[t]]$.

Any morphism of functor of groups $\, \phi \colon K \to A\,$ is the additive extension, $\phi = F_+$, of a unique series $\,F \in A(\pt)[[t]]$.
\end{teo}

\demo Again, let us only detail the case of a morphism $\phi \colon K \to A\,$ of functor of groups. Due to Axiom 6,  
the values of $\,\phi \,$ on the tautological line bundles $\, \xi_d \to \mathbb{P}^d\,$ are certain polynomials:
$$ 
\phi (\xi_d) \in A (\mathbb{P}^d) = \bigoplus_{i=0}^d A(\pt)\, x^i_d \ .  
$$

As these maps $\,\phi \,$ are compatible with the restrictions $i^* \colon \mathbb{P}^d \to \mathbb{P}^{d-1}$, the polynomials $\,\{ \phi (\xi_d)\}_{d \in \mathbb{N}}\,$ define a series $\, F \in A(\pt) [[t]]$. And we conclude $\,\phi  = F_+$, because both morphisms coincide on the tautological bundles (Lemma \ref{UnicidadMorfismoGrupos}).

\qed

\begin{ejem}\label{EjemSm}
The additive extensions of the monomials $\,t^m\,$ are morphisms of groups:
$$ s_m \colon K \to A \quad , \quad s_m (E) = \alpha_1^m + \ldots + \alpha_r^m \ , $$ which, due to the Newton's formulae, can be written in terms of the Chern classes:
\begin{align*}
s_0(E) &= \rk ( E) \\
s_1 (E) &= c_1^A (E) \\
s_2 (E) &= c_1^A (E)^2 - 2 c_2^A(E) \\ \ldots
\end{align*}
\end{ejem}

\begin{obser}\label{SerieInversa} If the series $\,F\,$ is invertible, then $\,F_\times\,$ extends to a natural transformation of groups $\, F_\times \colon K \longrightarrow (A^* , \cdot )$ that takes values in the multiplicative group $A^*$ of units of the cohomology.
Moreover, for any element $a \in K(X)$,
$$ F_\times (-a) = F_\times (a)^{-1} = (F^{-1})_\times (a) \ . $$
\end{obser}

\subsubsection*{Morphisms of rings $K \to A$}

\begin{teo}\label{MorfismosAnillos}
Any morphism of functor of rings $\,\varphi \colon K \to A\,$ is the additive extension, $\varphi = F_+$, of a unique series $\, F \in A(\pt) [[t]]\,$ satisfying
\begin{equation}\label{MorfismoFGL}
F( f^A (x,y) ) = F(x) \cdot F(y) \ .
\end{equation}
\end{teo}

\demo A morphism of functor of rings $\,\varphi \colon K \to A\,$ is the additive extension of a unique series $\,F \in A(\pt) [[t]]$ (Theorem \ref{MorfismosDeGrupos}).  By the Splitting Principle, the extension $F_+$ is a morphism of rings if and only if for any pair of line bundles, $L_x$ and $L_y$,
$$F_+ (L_x) \cdot F_+ (L_y) \, = \, F_+ (L_x \otimes L_y ) \ ; $$  
that is to say, if the series $\,F\,$ satisfies:
$$ F(x) \cdot F(y)=F (c_1^A (L_x \otimes L_y ) ) =F( f^A (x,y) )  .$$ 
\qed

\begin{ejem}\label{EjMorfismos}
\begin{enumerate}
\item
Any morphism of functor of rings $\,\varphi \colon K \to K\,$ is  the additive extension, $\varphi = F_+$, of a unique series $\, F \in \mathbb{Z} [[t]]\,$ satisfying
\begin{equation*}
F( x + y - xy ) = F(x) \cdot F(y) \ .
\end{equation*}For any integer $\,j$, the series $\, (1-t)^{j} $ satisfies this condition.

The {\it Adams operation} $\,\psi^j \colon K \to K\,$ is the additive extension of the series\footnote{The minus sign in the exponent is chosen in order to agree with the classical definition.} $\,(1-t)^{-j}$. These operations are morphisms of rings $\,\psi^j \colon K(X)  \longrightarrow K (X)\,$ that, for a non-zero $j$, act on line bundles as follows:
$$ \psi^j (L_x) = (1-x)^{-j} = \left( 1 - \left( 1 - L_x^{-1} \right) \right)^{-j}  = L_x^{j} \ , $$ where 
$\,x = c_1^K (L_x) = 1 - L_x^{-1}$. In particular,  $\,\psi^{-1}(L)=L^*$.


\item 
If $A$ follows the additive law and $A(\pt)$ is a $\mathbb{Q}$-algebra,  the exponential series $$\exp (t) := \sum_n \frac{t^n}{n!}\,  \in A(\pt) [[t]] $$
fulfills equation (\ref{MorfismoFGL}), and therefore its additive extension $\ch \colon K \to A$ is a morphism of rings, called the {\it Chern character}.

\item 
Let $A(\pt)$ be a ring where $1, 2, \ldots , p-1$ are units, and that contains a nilpotent element $\varepsilon$ such that $\varepsilon^p = 0$ (for example, $k [\varepsilon]/ (\varepsilon^p)$ with $k$ a field with characteristic zero or bigger or equal than $p$). If, moreover, $A$ follows the additive law, then the twisted exponential
$$ \exp_{\varepsilon} (t) := \exp ( \varepsilon \, t) = 1 + \varepsilon \, t + \ldots + \frac{\varepsilon^{p-1}}{(p-1)!}\, t^{p-1} \in A(\pt) [t] $$
also satisfies equation (\ref{MorfismoFGL}), and hence its additive extension $\ch_{\varepsilon} \colon K \to A$ is a morphism of rings, that we call {\it twisted Chern character}.

\end{enumerate}
\end{ejem}

\section{Integral formula for closed inmersions}

Let $\,A\,$ and $\,\bar{A}\,$ be oriented cohomology theories, and let us write 
$$ x_d := c_1^A( \xi_d)  \quad , \quad \bar{x}_d := c_1^{\bar{A}} (\xi_d) \  . $$

\begin{propo}\label{SerieAsociada} Let $\,\phi \colon A \to \bar{A}\,$ be a morphisms of functor of groups. There exists a unique series $\, \SS \in \bar{A}(\pt)[[t]]\,$, called the {\it series associated to} $\,\phi\,$, such that 
\begin{equation*}
 \phi (x_d ) \, = \, \SS (\bar{x}_d)\, \bar{x}_d  \ ,  
\end{equation*}
for all $\, d \in \mathbb{N} \cup \{0\}$. 
\end{propo}

\demo As before, recall that $A(\mathbb{P}^d)\, = A(\pt)[u] /(u^{d+1})$, and that the morphisms of groups $\,\phi \colon A(\mathbb{P}^d) \to \bar{A}(\mathbb{P}^d)\,$  are compatible with restrictions to the hyperplanes $i^* \colon A(\mathbb{P}^{d}) \to A(\mathbb{P}^{d-1})$.  

This implies the existence of a series $\, a_0 + a_1 t + \ldots \in \bar{A}(\pt) [[t]]\,$ such that 
$$
\phi (x_d) = a_0 + a_1\, \bar{x}_d + a_2\, \bar{x}_d^2 + \ldots \ .
$$

For $\,d=0$, we have $\,x_0 = 0 = \bar{x}_0\,$ and it follows $\, 0 = \phi (0) = a_0$. 

\qed

\begin{ejem}\label{EjemSeries}
\begin{enumerate}
\item If $A$ follows the additive law, then the series associated to $ s_m \colon K \to \bar{A}$ is the monomial $(-t)^{m-1}$:
$$ s_m (x_d) \, = \, s_m (1 - \xi^*_d) \, = \, 0 - \bar{y}_d^m \, = \,- (-\bar{x}_d)^m \, = \, (-\bar{x}_d)^{m-1} \bar{x}_d \ . $$

\item  The series associated to the Adams operation $\psi^j \colon K \to K $ is the polynomial
$$ \BB^j := \frac{1 - (1-t)^{j}}{t}\ = \ j - \binom{j}{2} t +  \ldots + (-1)^{j-1}t^{j-1} \ ,$$
because
$$\psi^j (x_d )  \, =\, \psi^j (1 - \xi_d^* ) \, =\, 1 - (\xi_d^*)^j = 1 - (1-x_d)^{j} \, = \, \left( \frac{1 - (1-x_d)^{j}}{x_d} \right) x_d  \ .
$$

\item If $A$ follows the additive  law and $A(\pt)$ is a $\mathbb{Q}$-algebra,  then
the series associated to the Chern character $\ch \colon K \to A$ is $ \frac{1-e^{-t}}{t}$:
$$
\ch (x_d )  \, =\, \ch (1 - \xi_d^* ) \, =\, 1 - e^{y_d} \, = \, 1 - e^{-x_d} \, = \, \left( \frac{1-e^{-x_d}}{x_d} \right)\, x_d \ ,
$$
(observe it is the inverse of the standard {\it Todd series} $ \Td (t) :=  \frac{t}{1-e^{-t}}$).

\item If $A$ follows the additive law and $A(\pt)$ is as in Example \ref{EjMorfismos}.3, then the series associated to the twisted Chern character $\ch_{\varepsilon} \colon K \to A$ is 
$$
\frac{1-\exp_\varepsilon (- t)}{t}= \varepsilon- \frac{\varepsilon^2}{2}\,t + \ldots + \frac{\varepsilon^{p-1}}{(p-1)!}\, (-t)^{p-2}  .
$$

\end{enumerate}
\end{ejem}






\begin{teo}\label{Main}
Let $\,\varphi \colon A \to \bar{A}\,$ be a morphism of functor of rings, with associated series $\, \SS \in \bar{A}(\pt)[[t]]\,$. For any closed inmersion $\, i \colon Z \hookrightarrow X\,$, the following square commutes:
\begin{equation*}\label{CuadradoEntero}
\vcenter{\xymatrix{
A(Z)  \ar[rrr]^-{i_*} \ar[d]_{  \varphi} & & & A(X) \ar[d]^{\varphi} \\
\bar{A} (Z) \ar[rrr]^-{\bar i_* ( \SS_{\times} (N_i) \cdot \raya ) } &  &  & \bar{A}(X),} }
\end{equation*}
where $\, N_i := i^* T_{X} - T_Z \,$ denotes the normal bundle of $i$. 

In other words, for any element $a \in A(Z)$,
\begin{equation}\label{EcuacionEntera}
\varphi \, i_*(a)  \, =\,  \bar i_*\left( \SS_{\times} (N_i) \cdot \varphi(a)\right) \ .
\end{equation}
\end{teo}

\demo The proof is established using the so called deformation to the normal bundle; that is to say, proving these two facts:

\begin{enumerate}
\item \emph{If equation (\ref{EcuacionEntera}) holds for the zero section $s\colon Z\to \tilde{N} :=\PP(1\oplus N_i)$ of the projective closure of the normal bundle, then it also holds for the closed immersion $i\colon Z\to X$}.

\item \emph{Equation (\ref{EcuacionEntera}) holds for the zero section $s \colon Z\to \tilde{E} :=\PP(1\oplus E)$ of the projective closure of any vector bundle $E\to Z$.}

\end{enumerate}

These statements are proved in closed analogy with points 1 and 2 of\cite[p. 332]{ExpositionesI}, essentially replacing $\bar i_*$,  $\bar s_{0*}$, and $\bar s_*$  in {\it loc. cit.} by  $\bar i_*\left( \SS_{\times} (N_i) \cdot \raya \right)$,   $\bar s_{0*}\left( \SS_{\times} (N_i) \cdot \raya \right) $ and $\bar s_*\left( \SS_{\times} (E) \cdot \raya \right) $,  respectively.
\qed

\begin{coro}[Adams-Riemann-Roch without denominators]\label{CoroRRAsindeno}
For  any closed inmersion $\, i \colon Z \hookrightarrow X\,$ and any $\,j \in \mathbb{Z}$, the square below commutes:
\begin{equation*}
\vcenter{\xymatrix{
K(Z)  \ar[rrr]^-{i_*} \ar[d]_{ \psi^j} & & & K(X) \ar[d]^{\psi^j} \\
K(Z) \ar[rrr]^-{i_* (\theta^j (N_i^*) \cdot \raya ) } & &  & K(X)  \ , }}
\end{equation*} where $N^*_i$ is the conormal bundle and $\theta^j \colon K^+ \to (K, \cdot ) $ is the morphism of functor of monoids defined on line bundles as
$$ 
\theta^j (L) := 1 + L + \ldots + L^{j-1} = \frac{1-L^j}{1-L} \ . 
$$ 

\end{coro}

\demo  Consider the multiplicative extension $ \BB^j_\times \colon K^{+} \to (K, \cdot ) $ of the polynomial $\frac{1 - (1-t)^{j}}{t}$ (Example \ref{EjemSeries}.2).
On any line bundle,
\begin{equation*}
\BB^j_\times (L_x) \, =\,  \frac{1 - (1-x)^j}{x} \, =\, \frac{1 - (L_x^{*})^{j}}{1- L_x^{*}} \,  \\= \,\theta^j (L_x^*) \ ,
\end{equation*} 
and the Splitting Principle then implies $\, \BB^j_\times ( N_i ) \, = \, \theta^j ( N_i^* )$. 

With this equality in mind, the thesis is just a particular case of Theorem \ref{Main}.
\qed



Applying Theorem \ref{Main} to the twisted Chern character of Example \ref{EjMorfismos}.3 immediately yields

\begin{coro}\label{FORMULA_NUEVA}
Let $A$ be a cohomology theory following the additive law, and such that  $A(\pt)$ is a ring where $1, 2, \ldots , p-1$ are units and that contains a nilpotent element $\varepsilon$ of order $p$ (\cf Example \ref{EjMorfismos}.3).  For  any closed inmersion $\, i \colon Z \hookrightarrow X\,$, the square below commutes:
\begin{equation*}
\vcenter{\xymatrix{
K(Z)  \ar[rrr]^-{i_*} \ar[d]_{\ch_{\varepsilon}} & & & K(X) \ar[d]^{\ch_{\varepsilon}} \\
A(Z) \ar[rrr]^-{i_* ( \TT_{\varepsilon} (N_i) \cdot \raya ) } & & & A(X)  \ , }}
\end{equation*} where $\TT_{\varepsilon} \colon K^+ \to A$ is the multiplicative extension of the (non-invertible) polyomial $\varepsilon - \frac{\varepsilon^2}{2}\,t +  \ldots + \frac{\varepsilon^{p-1}}{(p-1)!}\, (-t)^{p-2} $. 
\end{coro}

\section{Formula for projective morphisms}

The Riemann-Roch formula for arbitrary projective morphisms (Theorem \ref{RR}) will be a consequence of two facts. 

Firstly, an elegant statement due to I. Panin, whose proof encapsulates most of the arguments used in the standard proofs of the Grothendieck-Riemann-Roch theorem (\cf \cite[Sect. 2]{ExpositionesI},  \cite{Panin-Smirnov}):

\medskip
\noindent {\bf Panin's Lemma:}
{\it  A morphism of functors of rings $\,\varphi \colon A \to \bar{A}\,$ is a morphism of cohomologies if and only if its associated series is the identity.
}
\medskip

Observe for a moment that, if a cohomology $A$ does not follow the multiplicative group law, then there are no morphisms of cohomologies $K \to A$ (Corollary \ref{CorolarioSerie}).  Otherwise,

\begin{coro} If $A$ follows the multiplicative law,  then the additive extension of 
$\frac{1}{1-t} = 1 + t + t^2 + \ldots $ is the unique morphism of cohomology theories $\, K \to A$.
\end{coro}

\demo The additive extension of $F:= (1-t)^{-1}$ defines a natural transformation of rings $\,F_+\colon K \to A\,$ because it satisfies condition (\ref{MorfismoFGL}) of Theorem \ref{MorfismosAnillos}.   It is a morphism of cohomologies because its associated series is the identity (Example \ref{RelationXeY}.3):
\begin{equation*}
F_+ (x_d)\, =  \, 1 - F_+ (\xi_d^*) \, =\, 1 - F (y_d) \, = \, \frac{-c_1^A(\xi^*_d)}{1- c_1^A(\xi^*_d)} \, = \,  c_1^A(\xi_d) \ .
\end{equation*}

Uniqueness follows becuase $F_+$
preserves Chern classes,  so its value on the tautological bundles $\,\xi_d\,$ is completely determined
$$ F_+ ( \xi_d) = F_+ (1 - c_1^K (\xi_d^*) ) = 1 - F_+ ( y_d) = 1 - c_1^A (\xi_d^*) \ . $$ 
\qed

Secondly, the possibility of modifying direct images of cohomology theories:

\begin{defi} Let $\bar{A}$ be an oriented cohomology theory and let $\,\SS \in \bar{A}(\pt)[[t]]^*\,$ be an invertible series. The {\it modified direct image} $\,f_*^{\SS}\,$ of a projective morphism $f \colon Y \to X$ is defined as
$$f^{\SS}_{* } \, := \, \SS_{\times} (T_X) f_* \left( \SS_{\times} \left( T_Y \right)^{-1} \cdot \raya \right)\, =\, f_* \left( \SS_{{\times}} (-T_f) \cdot \raya \right)   \ ,$$ 
where $T_f := T_Y - f^* T_X \in K(Y)$ is the virtual relative tangent bundle. That is to say, $\, f^\SS\,$ is the dotted arrow that makes this square commutative:
$$
\xymatrix{
\bar{A}(Y) \ar[d]_{\SS_\times (T_Y)^{-1}\, \cdot  }  \ar@{..>}[rr]^-{f_*^{\SS} }  & & \bar{A}(X)  \\
\bar{A}(Y)  \ar[rr]^-{f_*}   & & \bar{A}(X) \ar[u]_{\cdot \, \SS_\times (T_X) } \ . }
$$
\end{defi}

Direct computation shows (\cite[Thm. 2.3]{ExpositionesI}):

\begin{propo}\label{Modificacion} The functor $\bar{A}$ endowed with  the modified direct images $f^\SS_*$ is an oriented cohomology theory, that we denote $(\bar{A} , f_{*}^{\SS})$.  Moreover, if we write $\bar{x}^{\SS}_d$ for the Chern class of $\xi_d$ in $(\bar{A} , f_{*}^{\SS})$, then
$$
\bar{x}^{\SS}_d \, = \, \SS(\bar{x}_d)\,\bar{x}_d \ . 
$$ 
\end{propo}

Now, the Riemann-Roch formula may be stated as follows:

\begin{teo}\label{RR} 
Let $\,\varphi \colon A\to \bar{A}\,$ be a morphism of functor of rings whose associated series $\SS \in \bar{A}(\pt)[[t]] $ is invertible. Then $\varphi$ is a morphism of cohomologies
$$ \varphi \colon \ A \ \longrightarrow \left( \bar{A} \, , \, f_{*}^{\SS} \right) \ .$$

In other words, for any projective morphism $\,f \colon Y \to X\,$ between smooth quasi-projective $k$-varieties, the diagram below commutes
$$
\xymatrix{
A(Y)  \ar[rrr]^-{f_*} \ar[d]_{\varphi} & & & A(X) \ar[d]^{\varphi} \\
\bar{A}(Y) \ar[rrr]^-{f_* \left( \SS_{\times} (-T_f) \cdot \raya \right) } & & & \bar{A} (X)}  \ .
$$
\end{teo}

\demo Due to Panin's Lemma, we only have to check that the series $\SS$ associated to $\varphi$ is the identity. But this  is a trivial consequence of Proposition \ref{Modificacion}: $$ \varphi ( x_d)\,  =\,  \SS (\bar{x_d}) \, \bar{x}_d\, =\, \bar{x}_d^{\SS} \ .$$ 
\qed

\begin{coro}[Adams-Riemann-Roch]\label{RRA} 
For any $j \in \mathbb{Z}$, the Adams operation $\psi^j$ is a morphism of cohomologies
$$ \psi^j \colon K \longrightarrow \left( K \otimes \mathbb{Z}[j^{-1}] \, , \,  f_{*}^{\BB^j} \right)\  . $$ 

Equivalently, for any projective morphism $\,f\colon Y\to X\,$ between smooth quasi-projective $k$-varieties, the following diagram commutes:
\begin{equation*}\label{RRAClasico}
\vcenter{\xymatrix{
K(Y)  \ar[rrr]^-{f_*} \ar[d]_{\psi^j} & & & K(X) \ar[d]^{\psi^j} \\
K(Y)\otimes \mathbb{Z}[j^{-1}]  \ar[rrr]^-{f_* \left( \theta^j (\Omega_f)^{-1} \cdot \raya \right)  } & & &K(X) \otimes \mathbb{Z}[j^{-1}]  \ , }}
\end{equation*}
where $\Omega_f = T_f^* \in K(Y) $ is the relative cotangent bundle.
\end{coro}

\demo The polynomial $\BB^j = j + \ldots$ is a unit in the ring of power series with coefficients in $K(\pt)[j^{-1}]$; hence, Theorem \ref{RR} readily implies the thesis.

The alternative description, using Bott's cannibalistic class $\theta^j$ and the cotangent bundle, is justified by the computation in the proof of Corollary \ref{CoroRRAsindeno}, that proves
$$\BB^j_\times ( -T_f ) \, = \, \theta^j ( -\Omega_f ) \, = \, \theta^j (\Omega_f )^{-1} \ . $$ \qed

\begin{coro}[Grothendieck-Riemann-Roch]\label{RRG} 
Let $A$ be an oriented cohomology following the additive group law and such that $A(\pt)$ is a $\mathbb{Q}$-algebra.  The Chern character is a morphism of  cohomologies
$$ \ch \colon K \longrightarrow \left( A  \, , \,  f_{*}^{\TT} \right) \ , $$ where $\TT = \frac{1-e^{-t}}{t}$ is its associated series (\cf Example \ref{EjemSeries}.3).

Equivalently, for any projective morphism $\,f\colon Y\to X\,$ between smooth quasi-projective $k$-varieties, the diagram below commutes:
$$
\xymatrix{
K(Y)  \ar[rrr]^-{f_*} \ar[d]_{\ch} & & &  K(X) \ar[d]^{\ch} \\
A(Y) \ar[rrr]^-{f_* \left( \Td (T_f) \cdot \raya \right)  } & & &  A(X)  \ , }
$$
where $\,\Td \colon K \to A \,$ is the multiplicative extension of the Todd series $\frac{t}{1-e^{-t}}$. 
\end{coro}

\demo The series $\TT = 1 - \ldots $ is a unit in $  \mathbb{Q} [[t]]$ , so the first statement is a direct consequence of Theorem \ref{RR}. 

The alternative (and classical) description using the Todd series is proved by the computation in Remark \ref{SerieInversa}, that implies
$$ \TT_\times (-T_f) \, = \,  (\TT^{-1})_\times (T_f)  \, = \, \Td ( T_f )  \ .  $$
\qed

\medskip
As for the last corollary, let us work with rational coefficients. Fix $\,j\in \mathbb{Z}$ and consider the following (invertible) series with rational coefficients:
$$ \BB^j = \frac{1 - (1-t)^{j}}{t}\quad , \quad \TT= \frac{1-e^{-t}}{t} \quad \mbox{ and } \quad  \SS^j:= \frac{1- e^{-jt}}{t}\, . $$




Let $CH_\QQ$ denote the Chow ring with rational coefficients, and consider the morphism of functor of rings: 
$$ \Phi^j :=\oplus_n \Phi^j_n \colon CH_\QQ \longrightarrow CH_\QQ \ ,  $$ where $\,\Phi^j_n\colon CH^n_\QQ \to CH^n_\QQ\,$ denotes  multiplication by $\,j^n\,$. 

\begin{coro}\label{TeoremaCuboMagico}
The following diagram is a commutative square of isomorphisms of cohomologies:
\begin{equation}\label{CuadradoCohomologias}
\vcenter{\xymatrix{K_\mathbb{Q}\ar[d]_-{\ch }\ar[rr]^-{\psi^j} & & \left( K_\mathbb{Q} , f^{\BB^j}_{* } \right) \ar[d]^-{\ch}\\
\left( CH_\mathbb{Q} , f^{\TT}_{*} \right) \ar[rr]^-{\Phi^j}& & \left( CH_\mathbb{Q}, f^{\SS^j}_{* } \right)}}
\end{equation} 

That is to say, the arrows above are isomorphisms of functor of rings such that, for any projective morphism $\,f\colon Y\to X\,$ between smooth, quasi-projective $k$-varieties, the cube below commutes:  
\begin{equation*}\label{CuboMagico}
\vcenter{\xymatrix@C-28
pt{ K(Y)_{\mathbb{Q}} \ar[dr]^{\Psi^j} \ar[dd]_{ \mathrm{ch} } \ar[rrr]^{ f_*} & & \qquad \phantom{A} \qquad  & K(X)_{\mathbb{Q}} \ar[dr]^{\Psi^j} \ar'[d][dd]^(.45){\mathrm{ch}}   &
\\
&  K(Y)_{\mathbb{Q}} \ar[dd]_(.30){\mathrm{ch}}\ar[rrr]^(.40){f_* \left( \theta^j (\Omega_f)^{-1} \cdot \raya \right) } & & & K(X)_{\mathbb{Q}}  \ar[dd]^{\mathrm{ch}}
\\
 CH(Y)_\mathbb{Q} \ar[dr]_{\Phi^j} \ar'[r][rrr]^(.40){f_* (\mathrm{Td} (T_f) \cdot \raya ) }& & & CH(X)_\mathbb{Q} \ar[dr]^{\Phi^j} &
\\
& CH(Y)_\mathbb{Q} \ar[rrr]^{f_* (\SS^j (-T_f) \cdot \raya ) } & & & CH(X)_\mathbb{Q} 
}} 
\end{equation*}
\end{coro}

\demo If we write $\,\mathrm{ch}_n :=\frac{s_n}{n!}$ (Example \ref{EjemSm}),  the commutativity of the square (\ref{CuadradoCohomologias}) is equivalent to the following equality of maps $\, K \to CH^n_{\QQ}$:
$$
\mathrm{ch}_n \circ \psi^j \, =\,  j^n \,  \mathrm{ch}_n \,\ .$$ 

Both terms are morphisms of groups, so it is enough to prove it for line bundles: 
\begin{align*}
\mathrm{ch}_n(\psi^j(L))&=\mathrm{ch}_n(L^{j}) = \frac{c_1 ^{CH}( L^{j})^n}{n!} = \frac{j^n\, c_1^{CH}(L)^n}{n!}=j^n\, \mathrm{ch}_n(L) \ . 
\end{align*}

Moreover,  the Chern character with rational coefficients is an isomorphism of rings (with inverse induced by $[Z]\mapsto [\O_Z]$, \cf \cite[XIV 4.2]{SGA6}), as so it is $\, \Phi^j\,$. 
Hence, $\psi^j$ is also an isomorphism of rings. 

Regarding the cube above,  Corollaries \ref{RRA}  and \ref{RRG} imply the commutativity of the top and back faces of the cube, respectively. 
For the front face,  the task is to prove that the Chern character is a morphism of cohomologies $\,\ch \colon \left(  K_\mathbb{Q} , f^{\BB^j}_{*} \right) \to \left( CH_\mathbb{Q} , f^{\SS^j}_{*} \right)\,$ and, for the bottom face,  that $\,\Phi^j \colon \left(  CH_\mathbb{Q} , f^{\TT}_{*} \right) \to \left( CH_\mathbb{Q} , f^{\SS^j}_{*} \right)\,$ is a morphism of cohomologies. 

Then, Panin's Lemma reduces these questions to check that their associated series are the identity:
\begin{align*}
\ch (x_d^{\BB^j}) &= \ch ( x_d \, \BB^j (x_d) ) = \ch(x_d)\frac{\ch (1 - (1-x_d)^j)}{\ch (x_d)} = 1 - \left( 1 - \ch (x_d) \right)^j \\
&= 1 - \left( 1 - (1- e^{-\bar{x}_d}) \right)^j 
=\, 1 - e ^{-j \bar{x}_d} \, = \, \SS^j (\bar{x}_d) \, \bar{x}_d \, = \, \bar{x}_d^{\SS^j} \ . \\
\Phi^j (\bar x_d^{\TT}) &= \Phi^j (\bar x_d \, \TT (\bar x_d) ) = \Phi^j(\bar x_d) \frac{1-\Phi^j (e^{-\bar x_d})}{\Phi^j(\bar x_d)} =1-e^{-j\bar x_d} = \, \bar{x}_d^{\SS^j} \ .
\end{align*}
\qed





\bibliography{GRR}{}

\providecommand{\bysame}{\leavevmode\hbox to3em{\hrulefill}\thinspace}
\providecommand{\MR}{\relax\ifhmode\unskip\space\fi MR }
\providecommand{\MRhref}[2]{%
  \href{http://www.ams.org/mathscinet-getitem?mr=#1}{#2}
}
\providecommand{\href}[2]{#2}
\begin{thebibliography}{10}

\bibitem{SGA6}
\emph{Th\'{e}orie des intersections et th\'{e}or\`eme de {R}iemann-{R}och},
  Lecture Notes in Mathematics, Vol. 225, Springer-Verlag, Berlin-New York,
  1971, S\'{e}minaire de G\'{e}om\'{e}trie Alg\'{e}brique du Bois-Marie
  1966--1967 (SGA 6), Dirig\'{e} par P. Berthelot, A. Grothendieck et L.
  Illusie. Avec la collaboration de D. Ferrand, J. P. Jouanolou, O. Jussila, S.
  Kleiman, M. Raynaud et J. P. Serre. \MR{0354655}

\bibitem{BorelSerre}
Armand Borel and Jean-Pierre Serre, \emph{Le th\'{e}or\`eme de
  {R}iemann-{R}och}, Bull. Soc. Math. France \textbf{86} (1958), 97--136.
  \MR{116022}

\bibitem{Duma}
Bertrand Duma, \emph{{Vers la forme g{\'e}n{\'e}rale du th{\'e}or{\`e}me de
  Grothendieck-{R}iemann-Roch}}, Theses, {Universit{\'e} Paris-Diderot - Paris
  VII}, September 2012.

\bibitem{Fulton}
William Fulton, \emph{Intersection theory}, second ed., Ergebnisse der
  Mathematik und ihrer Grenzgebiete. 3. Folge. A Series of Modern Surveys in
  Mathematics [Results in Mathematics and Related Areas. 3rd Series. A Series
  of Modern Surveys in Mathematics], vol.~2, Springer-Verlag, Berlin, 1998.
  \MR{1644323}

\bibitem{FultonLang}
William Fulton and Serge Lang, \emph{Riemann-{R}och algebra}, Grundlehren der
  mathematischen Wissenschaften [Fundamental Principles of Mathematical
  Sciences], vol. 277, Springer-Verlag, New York, 1985. \MR{801033}

\bibitem{Jouanolou}
J.~P. Jouanolou, \emph{Une suite exacte de {M}ayer-{V}ietoris en
  {$K$}-th\'{e}orie alg\'{e}brique}, Algebraic {$K$}-theory, {I}: {H}igher
  {$K$}-theories ({P}roc. {C}onf., {B}attelle {M}emorial {I}nst., {S}eattle,
  {W}ash., 1972), 1973, pp.~293--316. Lecture Notes in Math., Vol. 341.
  \MR{0409476}

\bibitem{LevineMorel}
M.~Levine and F.~Morel, \emph{Algebraic cobordism}, Springer Monographs in
  Mathematics, Springer, Berlin, 2007. \MR{2286826}

\bibitem{Manin}
Ju.~I. Manin, \emph{Lectures on the {$K$}-functor in algebraic geometry},
  Uspehi Mat. Nauk \textbf{24} (1969), no.~5 (149), 3--86. \MR{0265355}

\bibitem{ExpositionesI}
Alberto Navarro, \emph{On {G}rothendieck's {R}iemann-{R}och theorem}, Expo.
  Math. \textbf{35} (2017), no.~3, 326--342. \MR{3689905}

\bibitem{Panin-Smirnov}
I.~Panin, \emph{Riemann-{R}och theorems for oriented cohomology}, Axiomatic,
  enriched and motivic homotopy theory, NATO Sci. Ser. II Math. Phys. Chem.,
  vol. 131, Kluwer Acad. Publ., Dordrecht, 2004, pp.~261--333. \MR{2061857}

\bibitem{Panin}
Ivan Panin, \emph{{Oriented cohomology theories of algebraic varieties. II
  ({A}fter {I}. {P}anin and {A}. {S}mirnov)}}, Homology, Homotopy and
  Applications \textbf{11} (2009), no.~1, 349--405. \MR{2529164}

\bibitem{PinkRossler}
Richard Pink and Damian R\"{o}ssler, \emph{On the {A}dams-{R}iemann-{R}och
  theorem in positive characteristic}, Math. Z. \textbf{270} (2012), no.~3-4,
  1067--1076. \MR{2892938}

\bibitem{Riemann}
B.~Riemann, \emph{Theorie der {A}bel'schen {F}unctionen}, J. Reine Angew. Math.
  \textbf{54} (1857), 115--155. \MR{1579035}

\bibitem{Riou}
Jo\"{e}l Riou, \emph{Algebraic {$K$}-theory, {${\bf A}^1$}-homotopy and
  {R}iemann-{R}och theorems}, J. Topol. \textbf{3} (2010), no.~2, 229--264.
  \MR{2651359}

\bibitem{Roch}
G.~Roch, \emph{Ueber die {A}nzahl der willk\"{u}rlichen {C}onstanten in
  algebraischen {F}unctionen}, J. Reine Angew. Math. \textbf{64} (1865),
  372--376. \MR{1579304}

\bibitem{Smirnov}
A.~L. Smirnov, \emph{The {R}iemann-{R}och theorem for operations in the
  cohomology of algebraic varieties}, Algebra i Analiz \textbf{18} (2006),
  no.~5, 210--236. \MR{2301046}

\bibitem{Soule}
Christophe Soul\'{e}, \emph{Op\'{e}rations en {$K$}-th\'{e}orie
  alg\'{e}brique}, Canad. J. Math. \textbf{37} (1985), no.~3, 488--550.
  \MR{787114}

\bibitem{Weibel}
Charles~A. Weibel, \emph{The {$K$}-book}, Graduate Studies in Mathematics, vol.
  145, American Mathematical Society, Providence, RI, 2013, An introduction to
  algebraic $K$-theory. \MR{3076731}

\end{thebibliography}

\bibliographystyle{amsplain}

\end{document}